\documentclass{amsart}
\usepackage[colorlinks=true, citecolor=blue, linkcolor=red, urlcolor=blue]{hyperref}
\usepackage{graphicx}
\usepackage{amscd}
\usepackage{amsmath}
\usepackage{amsfonts}
\usepackage{cite}
\usepackage{amssymb}
\theoremstyle{plain}
\newtheorem{theorem}{Theorem}
\newtheorem{corollary}[theorem]{Corollary}
\newtheorem{lemma}[theorem]{Lemma}
\newtheorem{proposition}[theorem]{Proposition}
\theoremstyle{definition}
\newtheorem{example}[theorem]{Example}

\newtheorem{definition}[theorem]{Definition}

\newtheorem{remark}[theorem]{Remark}
\theoremstyle{remark}

\begin{document}
\title{$k-$Absorbing factorization in multiplicative lattices}

\author{Malik Tusif Ahmed and Tiberiu Dumitrescu}

\address{ Department of Mathematics and Statistics, Hazara University,
Mansehra, Pakistan}
\email{tusifahmed@hu.edu.pk, maliktusifahmed@outlook.com (Ahmed)}

\address{Facultatea de Matematica si Informatica, University of Bucharest, 14 A\-ca\-de\-mi\-ei Str., Bucharest, RO 010014, Romania}
\email{tiberiu\_dumitrescu2003@yahoo.com (Dumitrescu)}


\begin{abstract}
We study the multiplicative $C-$lattices generated by principal elements in which every proper element is a finite product of $k-$absorbing elements. We prove that such lattices have Krull dimension at most one. In the modular case, a finite direct product decomposition theorem is obtained.
\end{abstract}

\maketitle

\section{Introduction}
Let $R$ be a commutative unitary ring and $k\geq 2$ an integer.
Anderson and Badawi \cite{B,AB} introduced the concept of a (strongly) $k-$absorbing ideal  as a generalization of  prime ideals.
A proper ideal $Q$ of $R$ is said to be {\em strongly $k-$absorbing} (resp. {\em  $k-$absorbing}) if   whenever $Q$ contains a $(k+1)-$factor product 
$I_1\cdots I_{k+1}$ of ideals (resp. principal ideals) of $R$, it follows that $Q$ contains at least one $k-$factor subproduct of $I_1\cdots I_{k+1}$.
Recently, Secord \cite{S} proved that an ideal is $k-$absorbing iff it is strongly $k-$absorbing. The rings whose proper ideals are products of $k-$absorbing ideals were introduced and studied by the present authors, Muktar and Khadam  in \cite{MAD,ADK}.
Recently, Reinhart and Ulucak \cite{RU} studied the $2-$absorbing factorization property in the setup of Abstract Ideal Theory (i.e. in multiplicative lattices). 

In this  paper we 
study the $k-$absorbing  factorization property in multiplicative lattices. 
Let us recall briefly some basic terminology and facts from Abstract Ideal Theory (for further details, see \cite{A,AJ,AJ1,D,JJ}). 
A \emph{(multiplicative) lattice} is a complete lattice $(L,\leq)$, with least element $0$ and greatest element $1$, together with a commutative monoid   multiplication (with identity $1$) such that 
$ a\left(\bigvee_{\alpha}b_{\alpha}\right) =\bigvee_{\alpha}(ab_{\alpha})$  
 for every $a\in L$ and $\{b_{\alpha}\}_{\alpha}\subseteq L$. 
 As usual,  $\bigvee$  and $\bigwedge$ are the join and the meet in $L$ respectively.
 An element $c\in L$ is called \emph{compact} if, whenever $ c\leq\bigvee S $ for some subset $S\subseteq L$, there exists a finite subset $T\subseteq S$ such that $ c\leq\bigvee T. $ 
Denote by $L_\ast$ the set of compact elements of $L$. 
 A multiplicative lattice $L$ is called a \emph{$C-$lattice} if $1$ is compact, $L_\ast$ is closed under multiplication, and every element of $L$ is a join of compact elements. 
For $a,b\in L$, the residual quotient $(a:b)$ is defined by  $(a:b)=\bigvee\{c\in L\mid bc\leq a\}. $  Following \cite{D}, an element $e\in L$ is called: 
\emph{meet principal} if $a\wedge be=((a:e)\wedge b)e $ for all $a,b\in L$; 
\emph{join principal} if $(a\vee be:e)=(a:e)\vee b $ for all $a,b\in L$;
\emph{principal} if it is both meet principal and join principal. 
Every principal element of a $C-$lattice is compact (see \cite[Theorem 1.3]{A}). 
A $C-$lattice is said to be \emph{principally generated} if every element is a join of principal elements.

 {\em Throughout this paper by term {\em lattice} we mean   a $C-$lattice generated by principal elements. }     
Let $L$ be a lattice. 
An element $a\in L$ is \emph{proper} if $a<1$. A proper element $m$ is called \emph{maximal} if $ m<a\leq 1 $ implies  $a=1$. If $1$ is compact, then every proper element of $L$ is $\leq $ some maximal element. 
A lattice is said to be   \emph{local} if it has a unique maximal element. 
A proper element $p\in L$ is called \emph{prime} if $ ab\leq p $ (with $a,b\in L$) implies   $a\leq p$ or $b\leq p$. 
Every maximal element is prime. 
A lattice $L$  is a \emph{domain} if $0$ is a prime element. 
For a proper element $x\in L$, let $\operatorname{Min}(x)$
be the set of {\em minimal prime elements} above $x$.
The \emph{radical} $\sqrt{x}$ of an element $x\in L$ is defined by  
$\bigvee\{c\in L_\ast\mid c^n\leq x$ for some $n\geq 1 \}.$
We have $\sqrt{x}=\bigwedge Min(x)$.      
A proper element $q\in L$ is called \emph{primary} if $(q:x)\leq q$ whenever $x\not\leq q$; in this case $p=\sqrt{q}$ is prime and $q$ is called {\em $p$-primary}.
For $r\in L$, the {\em factor lattice} $L/r$ is $\{b\in L\mid r\leq b\}$ whose  multiplication is $ c\circ d=(cd)\vee r $ for all $c,d\in L/r$. 

Let $S$ be a multiplicative subset of $L_\ast$. For each $x\in L$, define $x_S=\bigvee_{s\in S}(x:s)$. The  {\em localization of $L$ at $S$} is $L_S=\{x_S\mid x\in L\}$ which is a lattice under multiplication $a\circ b= (ab)_S$ and join $X\mapsto  (\bigvee X)_S$ for $a,b\in L$ and $X\subseteq L$. 
If $p$ is a prime element of $L$, then $S(p)=\{x\in L_\ast\mid x\not\leq p\} $ is a multiplicative   subset of $L_\ast$ and 
$L_{S(p)}$, $x_{S(p)}$ are denoted simply by $L_p$ and $x_p$ respectively. Then $L_p$  (called the {\em localization of $L$ at $p$}) is a local lattice.
We have $x=\bigwedge_{m\in Max(L)} x_m$ for each $x\in L$, where $Max(L)$ is the set of maximal elements of $L$.
 The (Krull) dimension $dim(L)$ of $L$ is the supremum of all   integers $n\geq 0$ for which there exists a chain of prime elements $p_0<p_1<\ldots <p_n$.
\\[3mm]

We summarize the results of our paper.
Let $L$ be a lattice. In Definition \ref{30}, we introduce the concepts of {\em $k-$absorbing element}, {\em $k-$absorbing factorization} and {\em $k-$absorbing factorization lattice} similar to the ring case. 
The $k-$absorbing element property may be checked on compacts (Proposition \ref{1}). The $k-$absorbing factorization lattice   property behaves well under factorization, localization and finite direct products (Propositions  \ref{2},  \ref{3},  \ref{18}).
In a   $k-$absorbing factorization lattice all elements have finitely many minimal primes and every element $x$ is $\geq $  some power of its radical, so $x$ is $\geq$ some product of elements in $Min(x)$
(Propositions  \ref{4},  \ref{5}).
While a $p$-primary element $x$ with $p^k \leq x$ is $k-$absorbing, the converse holds if $p$ is a comparable element (Proposition \ref{6}).
The main results of this paper are Theorems \ref{10}, \ref{12}, and \ref{17}.
Theorem \ref{10} asserts that 
for a $k-$absorbing factorization lattice $L$, the localization $L_m$ is a domain for each non-minimal maximal element $m$ of $L$; in particular, $L$ is at most one-dimensional.
From this fact, we derive that any $k-$absorbing factorization modular lattice $L$  is a finite direct product of one-dimensional $k-$absorbing factorization lattice domains  and/or  zero-dimensional  local lattices  with  maximal element  nilpotent of index $k$ (Theorem \ref{12}).
Finally, we show that a lattice domain $L$ has $k-$absorbing factorization iff  $Min(x)$ is finite for each proper nonzero $x\in L$ and 
$L_m$  has $k-$absorbing factorization for each  $m\in Max(L)$
(Theorem \ref{17}).

\section{Results}

Recall that throughout this paper by term {\em lattice} we mean   a $C-$lattice generated by principal elements. 
Let $a_1, \ldots , a_{n}$ be elements of a lattice $L$.
Let $m<n$. By an {\em $m$-subproduct} of $a_1a_2\cdots a_{n}$ we mean a product of the form $a_{i_1}\cdots a_{i_m}$ for some  
$1\leq i_1< i_2< \cdots <i_m\leq n$. In the sequel, 
{\em $k$ is a positive integer}.

\begin{definition}\label{30}
Let $L$ be a lattice and $x\in L$ be a proper  element.

$(i)$ We say that $x$ is   a {\em $k-$absorbing element}, if whenever $a_1\cdots a_{k+1}\leq x$ (with  $a_i\in L$),  at least one 
$k$-subproduct of $a_1\cdots a_{k+1}$ is  $\leq x$.

$(ii)$ By a {\em $k-$absorbing factorization} of $x$ we mean an equality $x=x_1\cdots x_{n}$   whose factors $x_i$ are $k-$absorbing elements.

$(iii)$ If each proper element  of $L$ has at least one 
 $k-$absorbing factorization, we say that $L$ 
 {\em has $k-$absorbing factorization} (or that $L$ is a {\em $k-$absorbing factorization lattice}).
\end{definition}

Let $L(R)$ be the lattice ideal of a commutative ring $R$. 
Due to the main result of \cite{S}, an ideal $I$ of $R$ is a $k-$absorbing element of $L(R)$ iff $I$ is a $k-$absorbing ideal of $R$ in the terminology of \cite{AB}. So the rings $R$ with $L(R)$ a $k-$absorbing factorization lattice are exactly the FAF rings  of \cite{ADK}.
Next we show that  the $k-$absorbing property may be checked on compacts.

\begin{proposition}\label{1}
A proper  element  $x$ of  a lattice $L$ is  $k-$absorbing  iff whenever $c_1\cdots c_{k+1}\leq x$ with any $c_i$ compact in $L$,  
at least one $k$-subproduct of $c_1\cdots c_{k+1}$ is  $\leq x$.
\end{proposition}

\begin{proof} 
$(\Leftarrow)$ 
Deny, so suppose there exist $a_1,\ldots,a_{k+1}\in L$ such that  
product \\ $a_1a_2\cdots a_{k+1}\leq x$ but no $k$-subproduct  
 is $\leq x$. Since $$a_1a_2\cdots a_{k+1}=\bigvee \{ d_1d_2\cdots d_{k+1}|\ d_i \mbox{ compact } \leq a_i,\ i =1,\ldots,k+1\}$$
 and there are finitely many $k$-subproducts of $a_1a_2\cdots a_{k+1}$,
there exist some compacts $ c_i  \leq a_i, i=1,\ldots,k+1$ such that 
no $k$-subproduct of $c_1c_2\cdots c_{k+1}$   is $\leq x$, a contradiction.
The converse is trivial.
\end{proof}

\begin{proposition}\label{2}
If a lattice $L$ has $k-$absorbing factorization  and   $S\subseteq L_\ast$ a multiplicative set, then $L_S$ has $k-$absorbing factorization.
\end{proposition}

\begin{proof}
It suffices to prove that if $x$ is   $k-$absorbing in $L$, then $x_S$ is   $k-$absorbing in $L_S$, provided $x_S\neq 1$.
 We use Proposition \ref{1}.  Suppose that 
$(a_1a_2\cdots {a_{k+1}})_S\leq x_S$ with each $a_i$  compact in $L$.
By compactness, $(sa_1)a_2\cdots {a_{k+1}}\leq x$ for some $s\in S$; changing $a_1$ by $sa_1$, we may assume that $s=1$. Since $x$ is   $k-$absorbing in $L$, some $k$-subproduct of $a_1a_2\cdots {a_{k+1}}$ is $\leq x\leq x_S$.
\end{proof}

\begin{proposition}\label{3}
If a lattice $L$ has $k-$absorbing factorization  and   $1\neq r\in L$, then $L/r$ has $k-$absorbing factorization.
\end{proposition}

\begin{proof}
It is enough to prove the following fact: any $x\in L/r$ which is   $k-$absorbing in $L$ is  also $k-$absorbing in $L/r$.
Indeed, suppose that $a_1a_2   \cdots   a_{k+1}\vee r \leq x$ with 
$a_i\in L/r$.
As $x$ is $k-$absorbing in $L$, we get 
 $t\vee r\leq x$ for some $k$-subproduct $t$ of $a_1a_2\cdots a_{k+1}$.
\end{proof}

\begin{proposition}\label{18}
Two lattices $K$, $L$ have $k-$absorbing factorization iff their direct product $K\times L$ has $k-$absorbing factorization.
\end{proposition}
\begin{proof}
Implication $(\Leftarrow)$  follows from Proposition \ref{3}.
Conversely, it suffices to see that if $x$ and $y$  are $k-$absorbing elements of $K$ and $L$ respectively, the $(x,1)$ and $(1,y)$ are $k-$absorbing elements of $K\times L$.
\end{proof}

If $x$ is a proper element of a lattice $L$, denote by $Min(x)$ the set of minimal prime elements over $x$.

\begin{proposition}\label{4}
If $x$ is a proper element of a lattice $L$ and  $x_1\cdots x_n$ is a $k-$absorbing factorization of $x$, then $Min(x)$ has at most $nk$ elements.
\end{proposition}

\begin{proof}
Since $Min(x_1\cdots x_n)\subseteq \cup_{i=1}^n Min(x_i)$, it suffices to prove that $|Min(x)|\leq k$ when $x$ is   $k-$absorbing.
Deny, so suppose there exist $k+1$ distinct primes $p_1,\ldots , p_{k+1}\in Min(x)$. Pick compact elements $c_i\leq p_i$, $1\leq i \leq k$, such that $c_i\not\leq p_j$ for for all $j\in\{1,\ldots,k+1\}-\{i\}$. 
Next $p_i$ is minimal over $x$, so, changing $c_i$ by one of its powers, we may assume that   $s_ic_i\leq x$ for some compact $s_i\not\leq p_i$, $i=1,\ldots,k$.
Set $s:=s_1\vee \cdots \vee s_k$ to obtain $sc_1\cdots c_k\leq x$.
Note that $c_1\cdots c_k\not\leq x$ because $c_1\cdots c_k\not\leq p_{k+1}$. Similarly, if $q$ is the product $sc_1\cdots c_k$ with factor $c_i$ deleted, then $q\not\leq x$ because $q\not\leq p_i$. Consequently $x$ is not $k-$absorbing, a contradiction.
\end{proof}

The following   result is similar to \cite[Lemma 3.3]{RU}.

\begin{proposition}\label{5}
Let $x$ be a proper element of a   lattice $L$.

$(i)$ If $x$ has a $k-$absorbing factorization with $n$ factors, then $(\sqrt{x})^{nk}\leq x$. 

$(ii)$ If $x$ is radical, then $x$ is  $k-$absorbing iff $x$ has at most 
$k$ minimal primes. Thus the radical of 
a $k-$absorbing element is $k-$absorbing.

$(iii)$ If $x$ has a $k-$absorbing factorization, then 
$p_1p_2\cdots p_s\leq x$ for some  $p_i\in Min(x)$ (not necessarily distinct).
\end{proposition}
\begin{proof}
$(i)$ Suppose first that $n=1$ (i.e. $x$ is  $k-$absorbing). If $c$ is a compact element $\leq \sqrt{x}$, then $c^n\leq x$ for some $n$, hence   $c^k\leq x$ because $x$ is  $k-$absorbing. So $(\sqrt{x})^k = \bigvee \{c^k \ |\ c $ compact $ \leq \sqrt{x}\} \leq x$.
For $n\geq 2$, let  $y_1\cdots y_n$ be a $k-$absorbing factorization of $x$. We have
$$ (\sqrt{x})^{nk} \leq ((\sqrt{y_1})^k\wedge \cdots \wedge (\sqrt{y_n})^k)^n \leq (y_1\wedge \cdots \wedge y_n)^n \leq y_1\cdots y_n =x.$$ 
$(ii)$ Suppose that $Min(x)=\{p_1,\ldots,p_n\}$ with $n\leq k$; so
$x=p_1\wedge \cdots\wedge p_n$.
Let $a_1,\ldots,a_{k+1}\in L$ with $a_1\cdots a_{k+1}\leq\ {x}$. Then every $p_i$ is $\geq$ some $a_j$ (with $j$ depending on $i$). As $n\leq k$, some $k$-subproduct of $a_1\cdots a_{k+1}$ is $\leq x$. The converse follows from Proposition \ref{4}.
$(iii)$ Write  $\sqrt{x}=p_1\wedge \cdots\wedge p_n$ where $Min(x)=\{p_1,\ldots,p_n\}$ (cf. Proposition \ref{4}). 
By part $(i)$,  $(p_1 \cdots p_n)^k\leq \sqrt{x}^k\leq x.$
\end{proof}

The next lemma will be used in Proposition \ref{6}.

\begin{lemma}\label{20}
Let $p$ be a prime element of a chain lattice $L$. If $a,b\in L$ are $p$-primary elements with $ab\neq 0$, then  $ab$ is also  $p$-primary.
\end{lemma}
\begin{proof}
We have $\sqrt{ab} = \sqrt{a}\wedge \sqrt{b} = p$.
Let $x\in L$ be a principal element with $x \not\leq  p$. As   $L$ is a chain lattice, we get  $a \leq x$, so $a=x(a:x)=ax$  as $a$ is $p-$primary. 
Since $x$ is principal,  $((ab):x) = ((ab)x:x) = ab \vee (0:x) = ab$.  
Indeed, if $(0:x)\not\leq ab$, then  $ab \leq (0:x)$, so $ab=abx \leq (0:x)x=0$, a contradiction. So $ab$ is $p$-primary.
\end{proof}

\begin{example}\label{11}
In lemma above, condition $ab\neq 0$ is essential. 
Indeed, consider the valuation domain  $A=\mathbb{Z}_{(2)}+X\mathbb{Q}[[X]]$ and the factor chained ring $B=A/XA$. Let $P$ be the image of $X\mathbb{Q}[[X]]$ in $B$. Then $P$ is prime and $P^2=0$ is not primary because 
$2(X/2)B=0$, $2B\not\subseteq  P$ and $(X/2)B\neq 0$. Note that $0$ is not $2$-absorbing in $B$ (e.g. notice that $2^2(X/4)B=0$).
\end{example}

Let $L$ be a lattice and $x,p\in L$ with $p$ prime.
Recall  that the $k$-th {\em symbolic power} of $p$   is 
 the $p$-primary element $(p^k)_p$. And   $x$  is called  {\em comparable} if $x$ is comparable (under $\leq$) to any element of $L$.
The next result extends  
Proposition 2.2, Corollary 2.4 and Theorem 2.7 in \cite{ADK}.

\begin{proposition}\label{6}    
Let $p$ be a prime element  of a lattice $L$ and $x\in L$ with radical $p$. 

$(i)$ If $x$ is $p$-primary and $p^k \leq x$, then $x$ is $k-$absorbing.
In particular, the $n$-th  symbolic powers of $p$  for $1\leq n\leq k$     are $k-$absorbing.

$(ii)$ Conversely, if  $x$ is $k-$absorbing and $p$ is comparable, then $x$ is $p$-primary.

$(iii)$ In a  chain lattice $L$, the nonzero  $k-$absorbing elements are those of  the form $p^n\neq 0$ with $p$  prime    and $n\leq k$.
\end{proposition}
\begin{proof}
$(i)$ Let $y_1\cdots y_{k+1}\leq x$ with all factors $y_i$ in $L$. If all $y_i$'s are $\leq p$, then $y_1\cdots y_{k}\leq p^k \leq x$; if not, say $y_{k+1}\not\leq p$, then $y_1\cdots y_{k}\leq x$ as $x$ is $p$-primary. 

$(ii)$ Suppose that  $y,z$ are principal elements with $yz\leq x$ and $z\not\leq p$. As $p$ is comparable we get $y\leq p< z^{k-1}$, so $y=z^{k-1}t$ for some  $t\in L$.  As $x$ is $k-$absorbing, $yz=z^kt\leq x$ and $z\not\leq p\geq x$,
 we get $y=z^{k-1}t\leq x$. So $x$ is $p$-primary. 
 
 $(iii)$ Suppose $x$ is a nonzero $k-$absorbing element of $L$. Since $L$ is a chain lattice, $\sqrt{x}=p$ is prime.
 By $(ii)$ and Proposition \ref{5}, $x$ is $p$-primary and $p^n\leq x$ with $1\leq n\leq k$; choose   $n$ minimal with this property.
Pick a principal element $y\leq p^{n-1}$ with $y\not\leq x$. As $L$ is a chain lattice, we get $x=yh$ with $h\leq p$ (since $x$ is $p$-primary). Then  $x=yh\leq p^n$, so $x=p^n$. Conversely, an element of the form $p^n\neq 0$ with $p$  prime    and $n\leq k$
is $p$-primary (cf. Lemma \ref{20}), so it is $k-$absorbing cf. $(i)$.
\end{proof}

\begin{corollary}
A chain lattice domain has $k-$absorbing factorization iff it is a Dedekind lattice (i.e. its proper elements are finite products of prime elements).
\end{corollary}

\begin{example}
Assertion $(ii)$ Proposition \ref{6} above fails if the lattice is not generated by principal elements. Indeed, let
$L$  be the chain lattice $0< a < b < c < 1$ with multiplication 
$c^2=c$, $a^2=b^2=ab=ac=bc=a$.
Then  $a$ is $2$-absorbing with prime radical $b$ but is not primary as $bc=a$.
\end{example}

We give an example of $2$-absorbing nonprimary element.

\begin{example}\label{e1}
Let $L$ be a lattice and  $p,m\in L$ such that $p$ is  prime and $m$ maximal. Suppose that $p<m$ and $p\not\leq m^2$. Then $x:=p\wedge m^2$ is 
$2$-absorbing but not primary.
Indeed, $x$ is not primary because
$\sqrt{x}=p$, $m^2\not\leq p$, $p\not\leq x$  and $m^2p\leq x$. Suppose that $abc\leq x$ with $a,b,c\in L$. We may assume that $a\leq p$ and 
$a\not\leq m^2$. As $m^2$ is $m$-primary, we get $bc\in m$, so say $b\leq m$, hence $ab\leq x$.
In particular, $(X^2,2X)$ is a $2$-absorbing nonprimary ideal of 
$\mathbb{Z}[X]$.
\end{example}

The following three theorems are the main results of this paper.

\begin{theorem}\label{10}
A $k-$absorbing factorization lattice $L$ is at most one-dimensional. Moreover, $L_m$ is a domain for each non-minimal maximal element $m\in L$.
\end{theorem}
\begin{proof}
By contrary, {\em suppose that $dim(L)\geq 2.$ }
As the hypothesis is stable under factorization and localization, we may assume that $L$ is a  local lattice domain with  maximal element  $m$ having
 a prime element $0< p < m$. 
Pick a principal element $y\leq m$ with $y\not\leq p$. 
Shrinking $m$, we may assume that $m=\sqrt{p\vee y}$.
 Suppose  $m=m^2$.
As $\sqrt{p\vee y^2}=m$, we get $m=p\vee y^2$ (cf. part $(b)$ of  Proposition \ref{5}). By our choice of $y$ we get $1=(m:y)=(p\vee y^2):y=(p:y)\vee y=p\vee y$, a contradiction. 
{\em Thus  $m\neq m^2$.}
 Now we prove 
 
 {\em Claim $1$}: 
 every prime element $q < m$   is $\leq m^k.$\\
Deny, so there exists a prime element $q<m$ such that $q\not\leq m^k$. 
Pick a principal element $z\leq m$ with $z\not\leq q$. 
Let $t_1\cdots t_l$ be a  $k-$absorbing factorization of 
$q\vee z^{k^2}$; as $q\not\leq m^k$, we have $l<k$.
Since $z^{k^2}\leq t_j$ (for each $j$), we get  $z^{k}\leq t_j$, so
$z^{kl}\leq t_1\cdots t_l = q\vee  z^{k^2}$.
By our choice of $z$  we get $1=(q\vee z^{k^2}):z^{kl}=(q:z^{kl})\vee z^{k^2-kl} = q\vee z^{k^2-kl}$, a contradiction. 
So Claim $1$ above  holds.

As $m\neq m^2$, there exists a  principal element $w\leq m$ with $w\not\leq m^2$.    
By Claim $1$  we have $\sqrt{w}=m$. 
As $L$ has $k-$absorbing factorization and $w\not\leq m^2$, $w$ is a $k-$absorbing element, so $m^k\leq w$, cf. Proposition \ref{6}.
Since the $p\leq m^k\leq w$ (see Claim $1$) and  $w \not\leq p$, we get $p=w^kp$. 
Pick a nonzero principal element $x\leq p$. As $x$ has  a $k-$absorbing factorization and $p$ is prime, it follows that $u\leq p$ for some   $k-$absorbing  factor $u$ of $x$.
From $u\leq p\leq w^k$, we get $u=w^kh$ for some  $h\in L$. 
Since $w^k\not\leq u$ (as $w\not\leq p$) and $u$ is   $k-$absorbing, we get $w^{k-1}h\leq u$, hence 
$u=w^kh=wu$, thus $x=wx$ because $u$ is a factor of $x$. 
As $x$ is principal, we get $1=(w x:x)=w\vee (0:x)$, which is a contradiction because $w\neq 1$ and $x\neq 0$. Since this contradiction is caused by our initial assumption  $dim(L)\geq 2$, it follows that 
$dim(L)\leq 1$.
To prove the ``moreover part", we may suppose that $L$ is local with maximal element $m_0$ and 
$p_0<m_0$ is a prime element. Our task is to show that necessarily $p_0=0$. But this can be done reiterating the proof above for $m=m_0$ and $p=p_0$.
\end{proof}

\begin{example}\label{e2}
Let   $L$ be the lattice  having underlying set $\{0,1,a,b,c,d \}$ ordered by $a\leq c\leq d$,  $b\leq c$ with multiplication  given by
$\wedge$. Then $a,b,c,d$ are prime elements, so each proper element is a product of prime elements and $dim(L)=2$. So Theorem \ref{10} fails if the lattice is not principally generated. 
\end{example}

\begin{corollary}\label{14}
In a $k-$absorbing factorization lattice $L$, 
the non-maximal prime elements are idempotent.
\end{corollary}
\begin{proof}
Let $p$ be non-maximal prime of $L$ and $m$ a maximal element.
Then $p_m=0$ if $p\leq m$ (cf. Theorem \ref{10}) and 
$p_m=1$ if $p\not\leq m$. So $p=p^2$.

\end{proof}
\begin{remark}\label{15}
A   local zero-dimensional lattice $(L,m)$   has $k-$absorbing factorization iff $m^k=0$, cf. Propositions \ref{5} and \ref{6}.
\end{remark}
 
 Recall that a lattice $L$ is  \emph{modular} if $ a\wedge(b\vee c)=b\vee(a\wedge c) $ whenever $a,b,c\in L$ and $b\leq a$.

\begin{theorem}\label{12}
Any $k-$absorbing factorization modular lattice $L$  is a finite direct product of 

$(i)$ one-dimensional $k-$absorbing factorization lattice domains  and/or

$(ii)$ zero-dimensional  local lattices  with  maximal element  nilpotent of index $k$.
\end{theorem}
\begin{proof}
By Proposition \ref{5} $(iii)$ applied for $x=0$, we get 
$p_1^{c_1}\cdots p_n^{c_n}=0$ for some distinct minimal primes $p_1,\ldots,p_n$. By Theorem \ref{10} the  $p_i$'s are pairwise comaximal. So by Chinese Reminder Theorem it follows that $L$ is isomorphic to the direct product 
$L/p_1^{c_1} \times  \cdots \times  L/p_n^{c_n}$. When $p_i$ is non-maximal, $L/p_i^{c_i}$ is a domain by Corollary \ref{14}. Use Remark \ref{15}.
\end{proof}

\begin{remark}\label{16}
Let $m,x$ be elements of a lattice $L$ such that $m$ is maximal and $x$ is $m$-primary; clearly $x_m=x$.
Start with a $k-$absorbing factorization $x=(q_1\cdots q_n)_m$ of $x$ in $L_m$ where we may assume that $q_i=(q_i)_m$ for each $i$. Then  $x\vee m^k\leq q_i\leq m$, so $q_1\cdots q_n$ is $m$-primary, thus $x=q_1\cdots q_n$ is a $k-$absorbing factorization of $x$ in $L$.
\end{remark}

\begin{theorem}\label{17}
A lattice domain $L$ has $k-$absorbing factorization iff 

$(i)$ $Min(x)$ is finite for each proper nonzero $x\in L$ and 

$(ii)$ $L_m$  has $k-$absorbing factorization for each  $m\in Max(L)$.
\end{theorem}
\begin{proof}
Direction $(\Rightarrow) $ follows from Propositions \ref{2} and \ref{4}.
$(\Leftarrow)$ By Theorem \ref{10}, we may assume that $L$ is 
one-dimensional. Let $x$ be a proper nonzero element of $L$ and $Min(x)=\{p_1,\ldots,p_n\}$; these are exactly the maximal elements $\geq x$. Then 
$x= x_{p_1}\wedge \cdots \wedge x_{p_n}
=x_{p_1}  \cdots   x_{p_n}$ and  each $x_{p_i}$ is $m_i$-primary.
Apply Remark \ref{16} to conclude.
\end{proof}

We end this paper with two examples.

\begin{example}  
Let $(L,m)$ be a local one-dimensional lattice domain whose principal elements form an atomic monoid with finitely many atoms $a_1,\ldots , a_n$.
Then $m=a_1\vee\cdots \vee a_n$, so $m$ is compact; so each nonzero element is $\geq$ some power of $m$. We may assume that 
 $m^{k_i}\leq a_i$  $(1\leq i\leq n) $ for some integers $2\leq k_1< \ldots <k_n$.
We claim that {\em $L$ has  $(k_1-1)k_n-$absorbing factorization}. Note that 
$\bigwedge_{s\geq 1} a^s=0$ for each atom $a$, otherwise $a^t\leq \bigwedge_{s\geq 1} a^s$ for some $t\geq 1$, so $a^t=a^{t+1}$, a contradiction. Let $0\neq x\leq m$ be an element of $L$. If $x\not \leq m^{k_1}$, it follows that some product of $k_1-1$ atoms is $\leq x$, 
then $m^{(k_1-1)k_n}\leq  x$, so $x$ is $(k_1-1)k_n-$absorbing. 
Suppose that $x  \leq m^{k_1}\leq a_1$. We get   $x=a_1x_1$ for some $x_1\in L$ and, iterating this process, we end up with an equality 
$x=a_1^tz$ with $z\not\leq m^{k_1}$ (because $\bigwedge_{s\geq 1} a_1^s=0$). As argued above, $x=a_1^tz$ is a $(k_1-1)k_n-$absorbing factorization of $x$. Thus $L$ has $(k_1-1)k_n-$absorbing factorization.
A specific example is the lattice ideal $L(A)$ of domain 
$A=\mathbb{Z}[2i]_{(2,2i)} $. Indeed here the  atoms  are $2A$, $2iA$, 
$(2\pm 2i)A$, the maximal element is $m=(2,2i)A$ and we have $m^2\leq 2A\wedge 2iA$, $m^3\leq (2\pm 2i)A$. So $L(A)$ has $3-$absorbing factorization.
\end{example} 

\begin{example}  
Let $(L,m)$ be a local one-dimensional lattice domain.
Assume that $L$ has a sequence of proper elements  $(a_n)_{n\geq 1}$  such that  $m^n\not\leq a_n\not\leq m^2$. 
Suppose that $a_n$ has a $k-$absorbing factorization. As $a_n\not\leq m^2$ and $m^n\not\leq a_n$, it follows that $a_n$   itself is $k-$absorbing with $k>n$.
So {\em $L$ has not $k-$absorbing factorization for any $k$}.
A specific example is the lattice ideal $L(A)$ of domain 
$A=\mathbb{Z}[\sqrt{-7}]_{(2,1+\sqrt{-7})}$. Indeed here the maximal element is $m=(2,1+\sqrt{-7})A$ and
$m^n\not\leq 2[(1+\sqrt{-7})/2]^nA\not\leq m^2$ for each $n$ (to see, extend the ideals in the integral closure of $A$).
\end{example}


\begin{thebibliography}{9999}

 \bibitem{ADK} M.~T. Ahmed, T.~Dumitrescu, and M.~A. Khadam,
\emph{Commutative rings with absorbing factorization},
Comm. Algebra \textbf{48} (2020), no.~12, 5067--5075.
 
\bibitem{A} D.~D. Anderson, \emph{Abstract commutative ideal theory without chain condition}, Algebra Univ. \textbf{6} (1976), 131--145.


 \bibitem{AB} D.~F. Anderson and A. Badawi, \emph{On $n$-absorbing ideals of commutative rings}, Comm. Algebra {\bf 39} (2011), 1646--1672.

 \bibitem{AJ} D.~D. Anderson and C. Jayaram, \emph{Principal element lattices}, Czech. Math. J. \textbf{46} (1996), 99--109.

\bibitem{AJ1} D.~D. Anderson and E.~G. Johnson, \emph{Dilworth's principal elements,} Algebra Univ. \textbf{36} (1996), 392--404.


\bibitem{B} A. Badawi, \emph{On 2-absorbing ideals of commutative rings}, Bull. Austral. Math. Soc. {\bf 75} (2007), 417-–429. 


\bibitem{D} R.~P. Dilworth, \emph{Abstract commutative ideal theory} Pac. J. Math. {1\bf 2} (1962), 481-–498.



\bibitem{JJ} C. Jayaram and E.~G. Johnson, \emph{S-prime elements in multiplicative lattices}, Period. Math. Hung. {\bf 31} (1995), 201-–208. 

\bibitem{MAD} M. Mukhtar, M.~T. Ahmed and T. Dumitrescu, \emph{Commutative rings with two-absorbing factorization}, Comm. Algebra {\bf 46} (2018), 970--978.

\bibitem{RU} A. Reinhart and G. Ulucak, \emph{Multiplicative lattices with absorbing factorization}, Commun. Algebra, {\bf 53} (2025), no.~12, 5101--5119.

\bibitem{S} S. Secord, \emph{Collections of an ideal: any n-absorbing ideal is strongly n-absorbing}, 
Comm. Algebra, {\bf 52} (2024), 2649--2669.
\end{thebibliography}
\end{document}